\documentstyle[11pt,newlfont,amscd,amssymb]{amsart}
\newcommand{\nc}{\newcommand}

\renewcommand{\coprod}{\sqcup}

\nc{\one}{\mbox{\bf 1}}
\nc{\invtensor}{\underset{\leftarrow}{\otimes}}
\nc{\rlarrows}{\begin{picture}(1,0.4)
                \put(0,-0.1){\makebox(1,0.2){$\leftarrow$}}
                \put(0,0.1){\makebox(1,0.2){$\to$}}
              \end{picture}}
\nc{\rra}{\begin{picture}(1,0.4)
                \put(0,-0.1){\makebox(1,0.2){$\lra$}}
                \put(0,0.1){\makebox(1,0.2){$\lra$}}
              \end{picture}}
\nc{\Left}{\mathbf L}  % for derived
\nc{\Right}{\mathbf R} % functors
\nc{\gr}{\operatorname{gr}}
\nc{\Ho}{\operatorname{Ho}}
\nc{\alt}{\operatorname{alt}}
\nc{\Sym}{\operatorname{Sym}}
\nc{\sym}{\operatorname{sym}}
\nc{\id}{\operatorname{id}}
\nc{\Der}{\operatorname{Der}}
\nc{\im}{\operatorname{Im}}
\nc{\Ker}{\operatorname{Ker}}
\nc{\coker}{\operatorname{Coker}}
\nc{\Col}{\operatorname{Col}}
\nc{\ter}{\operatorname{ter}}
\nc{\intl}{\operatorname{int}}
\nc{\out}{\operatorname{out}}
\nc{\val}{\operatorname{val}}

\nc{\TN}{{\cal N}}
\nc{\Norm}{\operatorname{N}}
\nc{\Nor}{\operatorname{N}}
\nc{\Tor}{\operatorname{Tor}}
\nc{\res}{\operatorname{res}}
\nc{\Stab}{\operatorname{Stab}}
\nc{\Hom}{\operatorname{Hom}}
\nc{\chom}{\CH\!o\!m}
\nc{\uhom}{\CH\!o\!m}
\nc{\End}{\operatorname{End}}
\nc{\holim}{\operatorname{holim}}
\nc{\dirlim}{\underset{\rightarrow}{\lim}\,}
\nc{\invlim}{\underset{\leftarrow}{\lim}\,}
\nc{\CB}{\operatorname{\bf CB}}
\nc{\com}{\operatorname{co}}
\nc{\Tot}{\operatorname{Tot}}
\nc{\Th}{\operatorname{Th}}
\nc{\Cech}{\check{C}}
\nc{\Spec}{\operatorname{Spec}}
\nc{\Spf}{\operatorname{Spf}}
\nc{\MC}{\operatorname{MC}}
\nc{\U}{\operatorname{U}}
\nc{\Diff}{{\cal D}\mbox{\em iff}}
\nc{\Mor} {{\cal M}or}
\nc{\Ob}{\operatorname{Ob}}
\nc{\cone}{\widehat}
\nc{\Coder}{\operatorname{Coder}}
\nc{\pr}{\operatorname{pr}}
\nc{\diag}{\operatorname{diag}}

\nc{\HH}{$H_1$}

\nc{\CHo}{{\cal{H}\mbox{\it{o}}}}
\nc{\Mod}{{\mathtt{mod}}}       
\nc{\Modf}{{\mathtt{modf}}}       
\nc{\Modg}{{\mathtt{modg}}}       
\nc{\Ab}{{\mathtt {Ab}}}          
\nc{\Alg}{{\mathtt {Alg}}} 
\nc{\Hoalg}{{\mathtt {Hoalg}}} 
\nc{\Valg}{{\mathtt {Viral}}} 
\nc{\Algf}{{\mathtt {Algf}}} 
\nc{\Algg}{{\mathtt {Algg}}} 
\nc{\Coalg}{{\mathtt {Coalg}}} 
         
\nc{\dgc}{{\mathtt{dgc}}}
\nc{\dgca}{{\mathtt{dgca}}}
\nc{\dgcu}{{\mathtt{dgcu}}}
\nc{\dgcuf}{{\mathtt{dgcuf}}}
\nc{\dgcf}{{\mathtt{dgcf}}}
\nc{\dgcg}{{\mathtt{dgcg}}}
\nc{\dgcc}{{\mathtt{dgccc}}}

\nc{\dgl}{{\mathtt{dglie}}}
\nc{\dgla}{{\mathtt{dgla}}}
\nc{\dglf}{{\mathtt{dglf}}}
\nc{\dglg}{{\mathtt{dglg}}}
\nc{\dga}{{\mathtt{dga}}}
\nc{\art}{{\mathtt {art}}}
\nc{\dgar}{{\mathtt {dgart}^{\leq 0}}}
\nc{\simpl}{\Delta^{\op}\Ens}
\nc{\Coll}{{\mathtt{Coll}}}
\nc{\Colim}{{\mathtt{colim}}}
\nc{\Kan}{{\mathtt {Kan}}}

\nc{\Grp}{{\mathtt {Grp}}}
\nc{\Cat}{{\mathtt {Cat}}}
\nc{\Ens}{{\mathtt {Ens}}}
\nc{\op}{{\operatorname{op}}}
\nc{\Op}{{\mathtt{Op}}}

\nc{\Lie}{{\mathtt{LIE}}}
\nc{\Com}{{\mathtt{COM}}}
\nc{\Ass}{{\mathtt{ASS}}}
\nc{\pa}{\partial}

\nc{\CA}{\cal A}
\nc{\CDD}{\cal D}
\nc{\CE}{\cal E}
\nc{\CF}{\cal F}
\nc{\CG}{\cal G}
\nc{\CH}{\cal H}
\nc{\CI}{\cal I}
\nc{\CJ}{\cal J}
\nc{\CL}{\cal L}
\nc{\CM}{\cal M}
\nc{\CO}{\cal O}
\nc{\CP}{\cal P}
\nc{\CQ}{\cal Q}
\nc{\CR}{\cal R}
\nc{\CS}{\cal S}
\nc{\CT}{\cal T}
\nc{\CU}{\cal U}
\nc{\CW}{\cal W}
\nc{\CZ}{\cal Z}
\nc{\fa}{\frak a}
\nc{\fg}{\frak g}
\nc{\fk}{\frak k}
\nc{\fh}{\frak h}
\nc{\fm}{\frak m}
\nc{\fn}{\frak n}
\nc{\fA}{\frak A}
\nc{\fB}{\frak B}
\nc{\fC}{\frak C}
\nc{\fI}{\frak I}
\nc{\fS}{\frak S}

\nc{\nen}{\newenvironment}
\nc{\ol}{\overline}
\nc{\ul}{\underline}
\nc{\lra}{\longrightarrow}
\nc{\lla}{\longleftarrow}
\nc{\Lra}{\Longrightarrow}
\nc{\Lla}{\Longleftarrow}
\nc{\Llra}{\Longleftrightarrow}
\nc{\hra}{\hookrightarrow}
\nc{\iso}{\overset{\sim}{\lra}}

\nc{\Thm}[1]{Theorem~\ref{#1}}
\nc{\Prop}[1]{Proposition~\ref{#1}}
\nc{\Lem}[1]{Lemma~\ref{#1}}
\nc{\Cor}[1]{Corollary~\ref{#1}}
\nc{\Conj}[1]{Conjecture~\ref{#1}}
\nc{\Claim}[1]{Claim~\ref{#1}}
\nc{\Defn}[1]{Definition~\ref{#1}}
\nc{\Exa}[1]{Example~\ref{#1}}
\nc{\Rem}[1]{Remark~\ref{#1}}
\nc{\Note}[1]{Note~\ref{#1}}

\nen{thm}[1]{\label{#1}{\bf Theorem.\ } \em}{}
\nen{prop}[1]{\label{#1}{\bf Proposition.\ } \em}{}
\nen{lem}[1]{\label{#1}{\bf Lemma.\ } \em}{}
\nen{klem}[1]{\label{#1}{\bf Key Lemma.\ } \em}{}
\nen{cor}[1]{\label{#1}{\bf Corollary.\ } \em}{}
\nen{conj}[1]{\label{#1}{\bf Conjecture.\ } \em}{}
\nen{claim}[1]{\label{#1}{\bf Claim.\ } \em}{}

\nen{defn}[1]{\label{#1}{\bf Definition.\ } }{}
\nen{exa}[1]{\label{#1}{\bf Example.\ } }{}

\nen{rem}[1]{\label{#1}{\em Remark.\ } }{}
\nen{note}[1]{\label{#1}{\em Note.\ } }{}
\nen{exer}[1]{\label{#1}{\em Exercise.\ } }{}

\begin{document}

%\input{preamble}
%  top matter
\title[]{Virtual operad algebras and realization of homotopy types}
\author{Vladimir Hinich}
\address{Dept. of Mathematics, University of Haifa,
Mount Carmel, Haifa 31905 Israel}

%\thanks{}
\maketitle

%\tableofcontents
\section{Introduction}

\subsection{}
Let $k$ be a base commutative ring, $C(k)$ be the category of complexes
of $k$-modules. The category of operads $\Op(k)$ in $C(k)$ admits
a closed model category (CMC) structure with quasi-isomorphisms as
weak equivalences and surjective maps as fibrations (see~\cite{haha}, Sect.~6
and also Section~\ref{digest} below).

Let now $\CO$ be a cofibrant operad. The main result of this note 
(see~\Thm{main}) claims 
that the category of $\CO$-algebras admits as well a CMC structure
with  quasi-isomorphisms as weak equivalences and surjective maps as 
fibrations. This allows one, following the pattern of~\cite{haha}, 5.4,
to construct the homotopy category of {\em virtual} $\CO$-algebras
for any operad $\CO$ over $C(k)$ as the homotopy category of $\CP$-algebras
for a cofibrant resolution $\CP\to\CO$ of the operad $\CO$.

The main motivation of the note was to understand the following main result 
of Mandell's recent paper~\cite{man}.

\subsection{}
\label{man-th}

{\bf Theorem.}
{\em The singular cochain functor with coefficients in $\ol{\Bbb{F}}_p$
induces a contravariant equivalence from the homotopy category of
connected $p$-complete nilpotent spaces of finite $p$-type to a full
subcategory of the homotopy category of $E_{\infty}$ 
$\ol{\Bbb{F}}_p$-algebras.
}

In his approach, Mandell realizes the homotopy category of 
$E_{\infty}$-algebras as a localization of the category of algebras over a
``particular but unspecified'' operad $\CE$. One of major technical problems
was that the category of $\CE$-algebras did not seem to admit a CMC structure.

We suggest to choose $\CE$ to be a cofibrant resolution of the 
Eilenberg-Zilber operad. Then according to~\Thm{main}, the category of
$\CE$-algebras admits a CMC structure. This considerably simplifies the proof
of~\Thm{man-th}.

\subsection{Content of Sections}
The main body of the note (Sections 2 -- 4) can be considered as a
complement to~\cite{haha} where some general homology theory of operad
algebras is presented. 

In Section~2 we recall some results of~\cite{haha} we need in the sequel.
In Section~3 we prove the Main theorem~\ref{main}. In Section~4
we present, using~\Thm{main}, a construction of the homotopy category
$\Valg(\CO)$ of virtual $\CO$-algebras.

In Section~5 we review the proof of Mandell's theorem~\cite{man}, stressing
the simplifications due to our~\Thm{main}.

\subsection{Acknowledgement}This work was made during my stay at the 
Max-Planck Institut f\"ur Mathematik at Bonn. I express my gratitude
to the Institute for the hospitality. I am also grateful to P. Salvatore
for a useful discussion.

\section{Homotopical algebra of operads: a digest of~\cite{haha}}
\label{digest}

In this Section we recall some results from~\cite{haha} and give some
definitions we will be using in the sequel.

\subsection{Category of operads}
\label{digest-1}
Let $k$ be a commutative ring and let $C(k)$ denote the category
of complexes of $k$-modules. 

Recall (\cite{haha}, 6.1.1) that the category $\Op(k)$ of operads in $C(k)$ 
admits a closed model category (CMC) structure in which weak equivalences 
are componenwise quasi-isomorphisms and fibrations are componentwise 
surjective maps.

Cofibrations in $\Op(k)$ are retractions of {\em standard cofibrations};
a map $\CO\to\CO'$ is a standard cofibration if 
$\CO'=\dirlim_{s\in\Bbb{N}}\CO_s$ with $\CO_0=\CO$ and each $\CO_{s+1}$ 
is obtained from $\CO_s$ by adding a set of free generators $g_i$ with 
prescribed values of $d(g_i)\in\CO_s$.

\subsection{Algebras over an operad}
Let $\CO\in\Op(k)$.  

The category of $\CO$-algebras is denoted by $\Alg(\CO)$.
For $X\in C(k)$ we denote by $F(\CO,X)$ the free $\CO$-algebra
generated by $X$.
 
For any $d\in\Bbb{Z}$ denote by $W_d\in C(k)$ the contractible complex
$$ 0\to k\ = \ k\to 0$$
concentrated in degrees $d,d+1$.

\subsubsection{}
\begin{defn}{h1}
An operad $\CO\in\Op(k)$ is called \HH -operad if 
for any $A\in\Alg(\CO)$ the natural map
$$ A\to A\coprod F(\CO,W_d)$$
is a quasi-isomorphism.
\end{defn}

\subsubsection{}
\begin{prop}{2.2.1}(see~\cite{haha}, Thm.~2.2.1) 
Let $\CO$ be an \HH -operad. Then the category of $\CO$-algebras admits
a CMC structure with quasi-isomorphisms as weak equivalences and surjective
maps as fibrations.
\end{prop}

\subsection{Examples}
\subsubsection{} First of all, not all operads are \HH -operads. In fact,
let $k=\Bbb{F}_p$, $\CO=\Com$ (the operad of commutative algebras).
Then the symmetric algebra of $W_d$ fails to be contractible in degree $p$.

\subsubsection{}
\begin{prop}{4.1.1}(see~\cite{haha}, Thm.~4.1.1) Any $\Sigma$-split
operad (see~\cite{haha}, 4.2) is \HH -operad.
\end{prop}

In particular, all operads over $k\supseteq\Bbb{Q}$ are \HH -operads.
Also, all operads of form $\CT^{\Sigma}$ where $\CT$ is an asymmetric
operad, in particular, $\Ass$ (see~\cite{haha}, 4.2.5), are \HH -operads.

\subsubsection{}
The main result of this note claims that any  cofibrant operad
is an \HH -operad. 

\subsection{Base change and equivalence}

Let $f:\CO\to\CO'$ be a map of operads. Then a pair of adjoint functors
\begin{equation}
f^*:\Alg(\CO)\to\Alg(\CO'):f_*
\label{inv-dir-images}
\end{equation}
is defined in a standard way.

\subsubsection{}
\begin{prop}{4.6.4}(see~\cite{haha}, 4.6.4.)
Let $f:\CO\to\CO'$ be a map of \HH -operads.
The inverse and direct image functors~(\ref{inv-dir-images}) induce the adjoint
functors
\begin{equation}
\Left f^*:\Hoalg(\CO)\to\Hoalg(\CO'):\Right f_*=f_*
\label{Ho-inv-dir-images}
\end{equation}
between the corresponding homotopy categories.
\end{prop}

\subsubsection{}
\begin{defn}{strong-eq}A map $f:\CO\to\CO'$ of operads is called
{\em strong equivalence} if for each $d=(d_1,\ldots,d_n)\in\Bbb{N}^n$
the induced map
$$ \CO(|d|)\otimes_{\Sigma_d}k\to \CO'(|d|)\otimes_{\Sigma_d}k$$
is a quasi-isomorphism.

Here $|d|=\sum d_i$ and $\Sigma_d=\Sigma_{d_1}\times\ldots\times\Sigma_{d_n}
\subseteq\Sigma_{|d|}$.
\end{defn}

\subsubsection{}
\begin{prop}{eq}
Let $f:\CO\to\CO'$ be a strong equivalence of \HH-operads. Then the functors
$\Left f^*,\ f_*$ are equivalences.
\end{prop}

In Section~\ref{appl} we will be using  the following  version of \Prop{eq}.

\subsubsection{}
\begin{prop}{pre-eq}
Let $f:\CO\to\CO'$ be a strong equivalence of operads. Suppose $\CO$ is
\HH -operad. Then for each cofibrant $\CO$-algebra $A$ the natural map
$$ A\to f_*(f^*(A))$$
is an equivalence.
\end{prop}

\subsubsection{}
\begin{rem}{}
A quasi-isomorphism of $\Sigma$-split operads
compatible with the $\Sigma$-splittings is necessarily a strong equivalence.
\end{rem}

Theorem~4.7.4 of~\cite{haha} actually proves \Prop{pre-eq}
and \Prop{eq} together with the last Remark.
\section{Main theorem}

\subsection{}
\begin{thm}{main}Any cofibrant operad $\CO\in\Op(k)$ is an \HH -operad.

In particular, the category of algebras $\Alg(\CO)$ over a cofibrant
operad $\CO$ admits a CMC structure with quasi-isomorphisms as weak 
equivalences and epimorphisms as fibrations.
\end{thm}

\subsection{Proof of the theorem}
\label{pf-main}

First of all, we can easily reduce the claim to the case
$\CO$ is standard cofibrant. In fact,
since $\CO$ is cofibrant, it is a retraction of a standard cofibrant
operad $\CO'$. Let 
$$\CO\overset{\alpha}{\lra}\CO'\overset{\pi}{\lra}\CO$$
be a retraction. Let $A$ be a $\CO$-algebra. We can consider $A$ as a 
$\CO'$-algebra via $\pi$. Then the map $A\to A\coprod F(\CO,M)$ is a 
retraction of the map $A\to A\coprod F(\CO',M)$. This reduces the theorem to 
the case $\CO$ is standard cofibrant.

\subsection{Standard cofibrant case}
Let $\CO=\dirlim_{s\in\Bbb{N}}\CO_s$ (see notation of~\ref{digest-1},
$\CO_0=0$) 
be a standard cofibrant operad. Let $\{g_i\}$, $i\in I$ be a set
of free (homogeneous) generators of $\CO$. 

Let a function $s:I\to\Bbb{N}$
be given so that $\CO_s$ is freely generated as a graded operad by 
$g_i$ with $s(i)\leq s$ and, of course, $dg_i\in\CO_{s(i)-1}$.

Let, finally, $\val:I\to\Bbb{N}$ and $d:I\to\Bbb{Z}$ be the valency and the 
degree functions defined by the condition $g_i\in\CO(\val(i))^{d(i)}$.

The collection $\CI=(I,s,\val,d)$ will be called {\em a type} of $\CO$.

Since we deal with free operads and free algebras, it is worthwhile to have 
an appropriate notion of tree. Fix a type $\CI=(I,s,\val,d)$.

Put $I^+=I\cup\{a,m\}$ ($a$ and $m$ will be special marks on some terminal
vertices of our trees) and extend the functions $\val:I\to\Bbb{N}$ and 
$d: I\to\Bbb{Z}$ to $I^+$ by setting $\val(a)=\val(m)=d(a)=d(m)=0$.

\subsubsection{}
\begin{defn}{tree} {\em A $\CI$-tree} is a finite connected directed graph 
such that any vertex has $\leq 1$ ingoing arrows; each vertex is marked by 
an element $i\in I^+$ so that $\val(i)$ equals the number of outgoing arrows
which are numbered by $1,\ldots,\val(i)$.
\end{defn}

The set of vertices of a tree $T$ will be denoted by $V(T)$.
Terminal vertices of a $\CI$-tree are the ones having no outgoing arrows.
In particular, all vertices marked by $a$ or by $m$ are terminal.

\subsubsection{}
\begin{defn}{propertree}
 A $\CI$-tree $T$ is called {\em proper} if the following
property (P) is satisfied.

(P) For any vertex $v$ of $T$  one of the possibilities (a)--(c) below
occurs:

(a) $v$ is terminal;

(b) $v$ admits an outgoing arrow to a non-terminal vertex;

(c)  $v$ admits an outgoing arrow to a vertex marked by $m$.
\end{defn}

We denote by $\CP(\CI)$ the set of isomorphism classes of proper $\CI$-trees.
The following obvious result justifies the notion of proper tree.

\subsubsection{}
\begin{prop}{}
Let  $\CO$ be a standard cofibrant operad of type $\CI=(I,s,\val,d)$,
A be a $\CO$-algebra and $M\in C(k)$. Then the coproduct $B:=A\coprod F(M)$
is given, {\em as a graded $k$-module}, by the formula
\begin{equation}
B=\bigoplus_{T\in\CP(\CI)}A^{\otimes a(T)}\otimes M^{\otimes m(T)}[d(T)]
\end{equation}
where $a(T)$ (resp., $m(T)$) is the number of vertices of type $a$
(resp., of type $m$) in $T$ and $d(T)=\sum_{v\in V(T)}d(v).$
\end{prop}

\subsubsection{}
\label{lex}
Let $\CW$ be the set of maps $\Bbb{N}\to\Bbb{N}$ having finite support. 
Endow $\CW$ with the following lexicographic order. For $f,g\in\CW$ we will 
say that $f>g$ if there exists a $s\in\Bbb{N}$ such that $f(s)>g(s)$ and
$f(t)=g(t)$ for all $t>s$.

The set $\CW$ well-ordered.

Our next step is to define a filtration of $B=A\coprod F(M)$ indexed by
$\CW$.

\subsubsection{}
\begin{defn}{weight}
Let $T\in\CP(\CI)$. 
The weight of $T$, $w(T)\in\CW$ is the function $\Bbb{N}\to \Bbb{N}$ which
assigns to any $s\in\Bbb{N} $ the number of vertices $v$ of $T$ whose
mark  $i\in I$ satisfies $s(i)=s$.
\end{defn}

Now we are able to define a filtration on $B$.
\subsubsection{}
Let $A, M, B=A\coprod F(M)$ be as above. For each $f\in\CW$ define
$$ \CF_f(B)=
\bigoplus_{T: w(T)\leq f}A^{\otimes a(T)}\otimes M^{\otimes m(T)}[d(T)].$$

The homogeneous components of the associated graded complex are defined as
$$ \gr_f^{\CF}(B)=\CF_f(B)/\sum_{g<f}\CF_g(B).$$

\subsubsection{}
\begin{prop}{mainprop}
1. For each $f\in\CW$ the graded submodule $\CF_f$ is a subcomplex of $B$.

2. One has $\CF_0=A$.

3. Suppose $M$ is a contractible complex. Then for each $f>0$ the
homogeneous components $\gr_f^{\CF}$ are contractible.
\end{prop}
\begin{pf}Obvious. 
\end{pf}

\subsubsection{}
\begin{cor}{cmc}
The natural map $A\to B=A\coprod F(\CO,M)$ is a quasi-isomorphism of 
complexes. This implies Main Theorem~\ref{main}.
\end{cor}
\begin{pf}Obvious.
\end{pf}

\section{Virtual algebras}

\subsection{}
\label{virtual-sugg}

\Thm{main} suggests the following definition.

Let $\CO\in\Op(k)$. The homotopy category of virtual $\CO$-algebras
$\Valg(\CO)$ is defined as $\Hoalg(\CP)$ where $\CP\to\CO$ is a 
cofibrant resolution of $\CO$ in the category of operads.

One should, however, do some work, to ensure the definition above
makes sense.

\subsection{Base change}

Any morphism $f:\CP\to\CQ$ of operads induces a pair of adjoint functors
\begin{equation}
f^*:\Alg(\CP)\rlarrows\Alg(\CQ):f_*.
\label{invdirim}
\end{equation}

\Thm{main} together with~\ref{4.6.4} give immediately the following

\subsubsection{}
\begin{prop}{4.6.4(b)}For any morphism $f:\CP\to\CQ$ of cofibrant operads
the adjoint functors~(\ref{invdirim})
induce a pair of adjoint functors
\begin{equation}
\Left f^*:\Hoalg(\CP)\rlarrows\Hoalg(\CQ):\Right f_*=f_*
\label{hoinvdirim}
\end{equation}
between the homotopy categories. 
\end{prop}

\subsubsection{}
\begin{prop}{comp}1. Let $f:\CP\to\CQ$ be a weak equivalence of cofibrant
operads. Then $f$ is a strong equivalence. In paticular,
the derived functors of inverse and direct 
image~(\ref{hoinvdirim}) establish an equivalence of the homotopy categories.

2. Let $f,g:\CP\to\CQ$ be homotopic maps between cofibrant operads.
Then there is an isomorphism of functors
$$ f_*,g_*:\Hoalg(\CQ)\to\Hoalg(\CP).$$
This isomorphism depends only on the homotopy class of the homotopy
connecting $f$ with $g$.
\end{prop} 

\begin{pf}1. Let  $d=(d_1,\ldots,d_n),\ |d|=\sum d_i$
and let $\Sigma_d=\prod\Sigma_{d_i}\subseteq\Sigma_{|d|}$. 

We have to check that the map
$$ \CP(|d|)\otimes_{\Sigma_d}k\to\CQ(|d|)\otimes_{\Sigma_d}k,$$
induced by $f$, is a quasi-isomorphism.

Since $\CP$ and $\CQ$ are cofibrant operads, $\CP(|d|)$ and $\CQ(|d|)$
are cofibrant as complexes of $k(\Sigma_{|d|})$-modules. Therefore, their
quasi-isomorphism is a homotopy equivalence of $k(\Sigma_{|d|})$-modules
and therefore is preserved after tensoring by $k$.

2. We present here a proof which is identical to the proof of 
Lemma 5.4.3(2) of~\cite{haha}.

Let $\CQ\overset{\alpha}{\lra}\CQ^I\overset{p_0,p_1}{\rra}\CQ$  be a path 
diagram for $\CQ$ (see~\cite{q}, ch.~1) so that $\alpha$ is an acyclic 
cofibration. 
Since the functors $p_{0*}$
and $p_{1*}$ are both quasi-inverse to an equivalence 
$\alpha_*:\Hoalg(\CQ^I)\to\Hoalg(\CQ)$, they are naturally 
isomorphic. 
Therefore, any homotopy $F:\CP\to \CQ^I$ between $f$ and $g$ defines an
isomorphism $\theta_F$ between $f_*$ and $g_*$. Let now
$F_0,F_1:\CP\to\CQ^I$ be homotopic. The homotopy can be realized by a map
$h:\CP\to\CR$ where $\CR$ is taken from a path diagram
\begin{equation}
\CQ^I\overset{\beta}{\lra}\CR\overset{q_0\times q_1}{\lra} 
\CQ^I\times_{\CQ\times\CQ}\CQ^I
\label{2nd-path}
\end{equation}
where $\beta$ is an acyclic cofibration, $q_0\times q_1$ is a fibration,
$q_i\circ h=F_i, i=0,1.$ Passing to the corresponding homotopy categories
we get the functors $q_{i*}\circ p_{j*}: \Hoalg(\CQ)\to\Hoalg(\CR)$ which 
are quasi-inverse to 
$\alpha_*\circ\beta_*: \Hoalg(\CR)\to\Hoalg(\CQ)$. 
This implies that $\theta_{F_0}=\theta_{F_1}$.
\end{pf}

\subsection{Virtual operad algebras}

Our construction of the category of virtual $\CO$-algebras follows
the construction of virtual modules in~\cite{haha}, 5.4.

Let $\Op^c(k)$ denote the category of cofibrant operads in $C(k)$.
For each $\CP\in\Op^c(k)$ let $\Hoalg(\CP)$ be the homotopy category of 
$\CP$-algebras. These categories form a fibred category $\Hoalg$ over
$\Op^c(k)$, with the functors $\Right f_*=f_*$ playing the role of
``inverse image functors''.

Let $\CO\in\Op(k)$. Let   $\Op^c(k)/\CO$ be the category of maps $\CP\to\CO$
of operads with cofibrant $\CP$. The obvious functor 
$$c_{\CO}:\Op^c(k)/\CO\to\Op^c(k)$$
assigns the cofibrant operad $\CP$ to an arrow $\CP\to\CO$.

\subsubsection{}
\begin{defn}{virtual}The (homotopy) category $\Valg(\CO)$ of virtual 
$\CO$-algebras is the fibre of $\Hoalg$ at $c_{\CO}$. In other words, an
object of $\Valg(\CO)$ consists of a collection $A_a\in\Hoalg(\CP_a)$
for each $a:\CP_a\to\CO$ in $\Op^c(k)/\CP$ and of compatible collection 
of isomorphisms $\phi_f:A_a\to f_*(A_b)$ given for every $f:\CP_a\to\CP_b$ in 
$\Op^c(k)/\CO$.
\end{defn}

\subsubsection{}
\begin{cor}{descr}Let $\alpha:\CP\to\CO$ be a weak equivalence of operads with
cofibrant $\CP$. Then the obvious functor
$$ q_{\alpha}:\Valg(\CO)\to\Hoalg(\CP)$$
is an equivalence of categories.
\end{cor}
\begin{pf}We will construct a quasi-inverse functor 
$q^{\alpha}:\Hoalg(\CP)\to\Valg(\CO)$. For this choose for any map
$\beta:\CQ\to\CO$ a map $f_{\beta}:\CQ\to\CP$ making the corresponding
triangle homotopy commutative. Then, for any $A\in\Hoalg(\CP)$ we define
$q^{\alpha}(A)$ to be the collection of $f_{\beta*}(A)\in\Hoalg(\CQ).$
According to~\Prop{comp}, the definition does not depend on the choice
of $f_{\beta}'s$.
\end{pf}

The corollary means that the homotopy category of virtual $\CO$-algebras
is really the category of algebras over a cofibrant resolution of $\CO$.

\subsubsection{}
Any map $f:\CO\to\CO'$ defines an obvious functor
$\Op^c(k)/\CO\to\Op^c(k)/\CO'$. This induces a direct image functor 
$$f_*:\Valg(\CO')\to\Valg(\CO).$$
According to~\Cor{descr}, this functor admits a left adjoint 
inverse image functor $f^*$ which can be calculated using cofibrant
resolutions for $\CO$ and $\CO'$. 

\subsection{Comparing $\Valg(\CO)$ with $\Hoalg(\CO)$}

\subsubsection{}
Suppose $k\supseteq\Bbb{Q}$. Let $\CO\in\Op(k)$ and let $f:\CP\to\CO$ 
be a cofibrant resolution of $\CO$. Both $\CO$ and $\CP$ admit a
$\Sigma$-splitting (see~\cite{haha}, 4.2.4 and 4.2.5.2.) Moreover,
the quasi-isomorphism $f$ preserves the $\Sigma$-splittings. Therefore,
the categories $\Valg(\CO)=\Hoalg(\CP)$ and $\Hoalg(\CO)$ are equivalent
by ~\cite{haha}, 4.7.4.

Thus, in the case $k\supseteq\Bbb{Q}$  virtual operad algebras give 
nothing new.

\subsubsection{}
Let $\CT$ be an ``asymmetric operad'' i.e. a collection of complexes
$\CT(n)\in C(k)$ (with no action of the symmetric group), associative 
multiplication
$$ \CT(n)\otimes\CT(m_1)\otimes\ldots\otimes\CT(m_n)\to\CT(\sum m_i)$$
and unit element $1\in\CT(1)$ satisfying the standard properties.

Let $\CO=\CT^{\Sigma}$ be the operad induced by $\CT$ (see~\cite{haha}, 4.2.1).

\begin{lem}{}Suppose $\CT(n)$ are cofibrant in $C(k)$ (for example,
$\CT(n)\in C^-(k)$ and consist of projective $k$-modules). 
Then the natural functor
$$ \Valg(\CO)=\Hoalg(\CP)\to\Hoalg(\CO)$$
induced by a(ny) resolution $\CP\to\CO$, is an equivalence of categories. 
\end{lem}
\begin{pf} It is enough to check
that the map $\CP(n)\to\CO(n)$ is a homotopy equivalence of 
$k(\Sigma_n)$-complexes for each $n$.

But $\CP(n)$ is cofibrant over $k(\Sigma_n)$ since $\CP$ is a cofibrant
operad; $\CO(n)=\CT(n)\otimes k(\Sigma_n)$ is cofibrant over $k(\Sigma_n)$ 
since $\CT(n)$ is cofibrant over $k$. This proves the claim.
\end{pf}

\subsubsection{}Although the categories $\Valg(\CO)$ and $\Hoalg(\CO)$ 
turn out to be equivalent in all examples of $\Sigma$-split 
operads we know (\cite{haha}, 4.2.5), we do not see any reason
why this should always be the case. No doubt, the category $\Valg(\CO)$
should always be used when it differs from $\Hoalg(\CO)$.

\section{Application: realization of homotopy $p$-types}
\label{appl}

Mandell's theorem ~\cite{man} on the realization of homotopy $p$-types
can be reformulated in terms of virtual commutative algebras. The advantage
of this approach is that we can work with the category of operad algebras
which has a CMC structure. This makes unnecessary a big part of~\cite{man}.

In this Section we review the proof Mandell's theorem~\ref{man-th}.

\subsection{Adjoint functors $C^*$ and $U$}
\subsubsection{}
Recall~\cite{hlha} that the cochain complex $C^*(X)$ of an arbitrary
simplicial set $X\in\simpl$ admits a canonical structure of algebra
over the Eilenberg-Zilber operad $\CZ$ which is weakly equivalent
to the operad $\Com$ of commutative algebras. Choose any cofibrant resolution 
$\CE$ of $\CZ$. The category of virtual commutative algebras $\Valg(\Com)$
is canonically equivalent to $\Hoalg(\CE)$.

\subsubsection{}For each commutative ring $k$ define
\begin{equation}
C^*(\_,k):(\simpl)^{\op}\to\Alg(k\otimes\CE)
\label{cochains}
\end{equation}
(here and below $\otimes$ means tensoring over $\Bbb{Z}$)
to be the functor of normalized $k$-valued cochains.

This functor admits an obvious left adjoint functor
\begin{equation}
U_k:\Alg(k\otimes\CE)\to(\simpl)^{\op}
\label{functor-u}
\end{equation}
given by the formula
\begin{equation}
U_k(A)_n=\Hom(A,C^*(\Delta^n,k))
\label{def-u}
\end{equation}

The pair of functors $C^*(\_,k)$
 and $U_k$ satisfies the requirements of
Quillen's theorem~\cite{q}, \S 4, Theorem 3.

Since the functor $C^*(\_,k)$ preserves weak equivalences, one therefore 
obtains a pair of derived adjoint functors
\begin{equation}
\Bbb{U}_k:\Valg(\Com)=\Hoalg(k\otimes\CE)\rlarrows\CHo: C^*(\_,k),
\label{adj-pair}
\end{equation}
$\CHo$ being the homotopy category of simplicial sets.

\subsection{}
Following ~\cite{man}, we call $X\in\simpl$ {\em $k$-resolvable} if the
natural map
$$ u_X: X\to\Bbb{U}_k C^*(X,k)$$
is a weak equivalence.

The following two lemmas allow one to construct resolvable spaces.

\subsubsection{}
\begin{lem}{thm-1.1}(\cite{man}, Thm. 1.1) Let $X$ be the limit of a 
tower of Kan fibrations 
$$\ldots\to X_n\to\ldots\to X_0.$$
 Assume that the canonical map from $H^*X$ to $\Colim\ H^*X_n$ is an 
isomorphism.  If each $X_n$ is $k$-resolvable, then $X$ is $k$-resolvable.
\end{lem}

\subsubsection{}
\begin{lem}{thm-1.2}(\cite{man}, Thm. 1.2)
Let $X,Y$ and $Z$ be connected
simplicial sets of finite type, and assume that $Z$ is simply connected.
Let $X\to Z$ and $Y\to Z$ be given, so that $Y\to Z$ is a Kan fibration. Then,
if $X,Y$ and $Z$ are $k$-resolvable  then so is the fibre product $X\times_ZY$.
\end{lem} 

\Lem{thm-1.1} follows form the fact that the functor $\Bbb{U}$ carries homotopy
colimits in $\Alg(\CE)$ into homotopy limits in $\simpl$. The proof
of \Lem{thm-1.2} is similar, but needs in addition ~\Prop{ascom} below
which can be also  easily deduced from~\Thm{main}.

Using the CMC structure on $\Op(k)$, one can embed the obvious map 
of operads $\Ass\to\Com$ into the following commutative diagram
\begin{center}
$$\begin{CD}
\Ass_{\infty}@>{\alpha}>> \CE \\
@VVV @V{\tau}VV \\
\Ass @>{\ol{\alpha}}>> \ol{\CE} @>{\pi}>> \Com
\end{CD}
$$
\end{center}
where $\Ass_{\infty}$ is the operad of $A_{\infty}$-algebras, $\alpha$
is a cofibration, $\pi$ is a weak equivalence and the square is 
cocartesian.

\subsubsection{}
\begin{prop}{ascom}(compare to~\cite{man}, Lemma 5.2). Let $A\to B$ and $A\to C$ be 
cofibrations of cofibrant $\CE$-algebras. Let $\ol{A}=\tau^*(A)$, 
and similarly for $\ol{B},\ol{C}$. Then the natural maps
$$ B\coprod^AC\overset{t}{\longrightarrow}\ol{B}\coprod^{\ol{A}}\ol{C}
\overset{r}{\longleftarrow}\ol{B}\otimes_{\ol{A}}\ol{C}$$
are quasi-isomorphisms in $C(k)$.
Here $t$ is induced by $\tau$ and $r$ is induced by the composition
$$  \ol{B}\otimes\ol{C}\to (\ol{B}\coprod^{\ol{A}}\ol{C})\otimes
(\ol{B}\coprod^{\ol{A}}\ol{C})\overset{\operatorname{mult.}}{\longrightarrow}
\ol{B}\coprod^{\ol{A}}\ol{C}.$$
\end{prop}
\begin{pf}
1. {\em $t$ is a quasi-isomorphism.} The functor $\tau^*$ commutes with 
colimits. Therefore, it is enough
to prove that the natural map $A\to\tau_*\tau^*(A)$ is a weak equivalence
for a cofibrant algebra $A$. According to~\ref{pre-eq},
it is enough to check that $\tau:\CE\to\ol{\CE}$ is a strong equivalence of 
operads. 

Since $\alpha$ is a cofibration, $\ol{\alpha}$ is a cofibration as well.
Therefore, both $\CE(n)$ and $\ol{\CE}(n)$ are cofibrant over 
$k\Sigma_n$.
Then the strong equivalence of $\CE$ and $\ol{\CE}$ follows from their 
weak equivalence.

2.  {\em $r$ is a quasi-isomorphism.}

Suppose $A$ is standard cofibrant and the maps $A\to B$, $A\to C$ are
standard cofibrations. Let $\{e_i,i\in I\}$, $\{e_j,j\in I\cup J\}$,
$\{e_k,k\in I\cup K\}$, be graded free bases of $A,B$ and $C$ respectively
(the index sets $I,J,K$ are disjoint).

The sets $I,J$ and $K$ are well-ordered and the differential of $e_i$
is expressed through $e_{i'}$ with $i'<i$.

Put $S=I\cup J\cup K$ with the order given by $i<j<k$ for 
$i\in I,j\in J,k\in K$. Let $\widetilde{S}$ be the set of maps $S\to\Bbb{N}$
with finite support and with the lexicographic order as in~\ref{lex}.

For $f\in\widetilde{S}$ denote $|f|=\sum_{s\in S} f(s)$.

The algebra $\ol{B}\coprod^{\ol{A}}\ol{C}$
has an obvious increasing filtration by subcomplexes $\{F_f\}$ 
indexed by $f\in\widetilde{S}$.

The homogeneous component of the associated graded complex for 
$f\in\widetilde{S}$ takes form
$$ \gr_f(F)=\ol{\CE}(|f|)\otimes_{\Sigma_f}e^f$$
where $e^f=\prod_{s\in S}e_s^{f(s)}$ and 
$\Sigma_f=\prod_{s\in S}\Sigma_{f(s)}$.

Define a filtration $\{F'_f\}$ of $\ol{B}\otimes_{\ol{A}}\ol{C}$ 
indexed by the same set $\widetilde{S}$. It is given by the formula
$$
F'_f=\bigoplus_{g<f}\ol{\CE}(|g|_1)\otimes\ol{\CE}(|g|_2)\otimes_{\Sigma_g}e^g
$$
where $|g|_1=\sum_{s\in I\cup J} g(s)$ and $|g|_2=\sum_{s\in K} g(s)$.
The homogeneous component  for $f\in\widetilde{S}$ is given by
$$ \gr_f(F')=\ol{\CE}(|f|_1)\otimes\ol{\CE}(|f|_2)\otimes_{\Sigma_f}e^f.$$

The map $r: \ol{B}\otimes_{\ol{A}}\ol{C}\to \ol{B}\coprod^{\ol{A}}\ol{C}$
is compatible with the filtrations. The corresponding map of the homogeneous
components
$$ \gr_f(r):\ol{\CE}(|f|)\otimes_{\Sigma_f}e^f\to
\ol{\CE}(|f|_1)\otimes\ol{\CE}(|f|_2)\otimes_{\Sigma_f}e^f
$$
is induced by the map 
\begin{equation}
\ol{\CE}(|f|_1)\otimes\ol{\CE}(|f|_2)\to\ol{\CE}(|f|)
\label{split-e}
\end{equation}
which is obviously quasi-isomorphism. The assertion then follows from
the observation that both the left and the right hand side of~(\ref{split-e})
are cofibrant over $k(\Sigma_f)$.
\end{pf}

\subsection{}
To construct $k$-resolvable spaces using~\ref{thm-1.1} and~\ref{thm-1.2}
one needs a space ``to start with''. This is the Eilenberg-Maclane
space $K(\Bbb{Z}/p,n)$. The key step in ~\cite{man} is the following
\subsubsection{}
\begin{thm}{base}(cf.~\cite{man}, Prop. A.7). The space  $K(\Bbb{Z}/p,n)$ is 
$k$-resolvable iff $k\supseteq\Bbb{F}_p$ and the frobenius $F:k\to k$
gives rise to a short exact sequence of abelian groups
\begin{equation}
0\to\Bbb{F}_p\to k\overset{1-F}{\lra}k\to 0.
\label{ses}
\end{equation}
\end{thm}
\begin{pf}
This is the most important part of Mandell's result and we cannot simplify
the original Mandell's proof. 

1. The main step is to construct an
explicit cofibrant resolution of $C:=C^*(K(\Bbb{Z}/p,n),\Bbb{F}_p)$ 
over $k:=\Bbb{F}_p$.

Let $k=\Bbb{F}_p$. Let $\CE$ be a cofibrant resolution of the operad $\Com$
over $\Bbb{Z}$.
Recall that $\CE$-algebra structure on $A$ gives rise to the action of
the generalized Steenrod algebra $\fB$ on $H(A\otimes\Bbb{F}_p)$ 
--- see~\cite{may}.

Let  $A$ be a chain complex of a topological space and let the operad $\CE$ be
endowed with a weak equivalence $\CE\to\CZ$ to the Eilenberg-Zilber operad
of~\cite{hlha}, so that $A$ becomes an $\CE$-algebra. Then the action 
of $fB$ on $H(A\otimes\Bbb{F}_p$ induces
an action of the (conventional) Steenrod algebra $\fA$ which is a quotient
of $\fB$ by the ideal generated by $P^0-1$, $P^0$ being the degree zero
generalized Steenrod operation.

Choose a fundamental cycle $e\in C^n$. This cycle defines
a map $\phi:\CE_{\Bbb{F}_p}\langle x\rangle\to C$ from the free 
${\Bbb{F}_p}\otimes\CE$-algebra with
a generator $x$ to $C$ sending $x$ to $e$. Since $P^0$ acts trivially
on $H(C)$, the cohomology class $P^0([x])-[x]$ of 
$\CE_{\Bbb{F}_p}\langle x\rangle$
(here $[x]$ is the cohomology class of $x$), belongs to the kernel
of $H(\phi)$. Choose a representative $z$ of the cohomology class 
$P^0([x])-[x]$ of $\CE_{\Bbb{F}_p}\langle x\rangle$. 

Finally, define $B=\CE_{\Bbb{F}_p}\langle x,y;dy=z\rangle$. This is the 
$\Bbb{F}_p\otimes\CE$-algebra obtained from the free algebra 
$\CE_{\Bbb{F}_p}\langle x\rangle$ by adding a variable to kill the cycle $z$
~\see{haha}, 2.2.2.

The map $\phi$ can be obviously extended to a map $\psi: B\to C$.

\begin{thm}{resolution}(see~\cite{man}, Thm. 6.2). The map $\psi$ is a 
quasi-isomorphism.
\end{thm}

The proof of the theorem given in~\cite{man}, Sect. 12, is based on a
study of free unstable modules over $\fB$ and $\fA$.

2. Once we have found a cofibrant resolution $B$ of the algebra $C$ of cochains
of $K(\Bbb{Z}/p,n)$, the life becomes very easy. 

We have to study the map $u_X:X\to\Bbb{U}_k(C^*(X,k))$ for $X=K(\Bbb{Z},n)$.

One has $ \Bbb{U}_k(C^*(X,k))=U_k(B_k)$ where $B_k=k\otimes_{\Bbb{F}_p}B=
\CE_k\langle x,y;dy=z\rangle$. Since the functor $U_k$ carries
cofibrations to Kan fibrations and colimits to limits, one has
a cartesian diagram of spaces
\begin{center}
$$
  \begin{CD}
    U_k(B_k)@>>>U_k(\CE_k\langle z,y; dy=z\rangle) \\
    @VVV   @VVV \\
    U_k(\CE_k\langle x\rangle)@>p>>
                    U_k(\CE_k\langle t\rangle)\\
  \end{CD}
$$
\end{center}

The vertical maps are Kan fibrations and $U_k(\CE_k\langle z,y; dy=z\rangle)$
is contractible since $\CE_k\langle z,y; dy=z\rangle$ is a 
contractible $k\otimes\CE$-algebra. 

Furthermore, $U_k(\CE_k\langle t\rangle)$ identifies easily with the
Eilenberg-Mac Lane space $K(k,n)$ and the map $p$ is induced
by $1-F:k\to k$ where $F:k\to k$ is the frobenius (\cite{man}, prop. 6.4,
6.5).

Then the long exact sequence of the homotopy groups for the fibration
$U_k(B_k)\to U_k(\CE_k\langle x\rangle)$ gives the long exact sequence
$$\ldots\to\pi_i(K(k,n+1))\to\pi_i(U_k(B_k))\to\pi_i(K(k,n))\overset{p}{\lra}
\pi_{i+1}(K(k,n+1))\to\ldots$$
where the map $p$ is induced by $1-F$.

Now, if the condition on $k$ is not fulfilled, $ U_k(B_k)$ is not an 
Eilenberg-Mac Lane space. If the sequence~(\ref{ses}) is exact, the natural
map $K(\Bbb{Z}/p,n)\to U_k(B_k)$ induces isomorphism of homotopy groups
and this proves the assertion.
\end{pf}

\end{document}